# A Note on Lipshitz's Lemma 3

Shaoshi Chen   and   Ziming Li[*][†]

Department of Mathematics, North Carolina State University, Raleigh
KLMM, Academy of Mathematics and Systems Science, Beijing

November 5, 2011

## Abstract

In this note, we give a remark on the proof of Lemma 3 by Lipshitz in [1]. This remark is motivated by the observation that the statement from line $-8$ to $-3$ on page 375 of [1] seems not completely correct.

## 1   An algebraic description of Lipshitz's Lemma

Let $K$ be a field of characteristic zero, and $K(x,y)$ be the field of rational functions in $x$ and $y$ over $K$. Denote by $\mathcal{R}_2$ the ring $K(x,y)\langle D_x, D_y \rangle$ of linear differential operators generated by $D_x$ and $D_y$ over $K(x,y)$, whose commutative rules are given by

$$D_x f = f D_x + \frac{\partial f}{\partial x} \quad \text{and} \quad D_y f = f D_y + \frac{\partial f}{\partial y} \quad \text{for all } f \in K(x,y).$$

Lemma 3 in [1] is an easy consequence of the following proposition.

**Proposition 1.1** *Let $I$ be a left ideal of $\mathcal{R}_2$. If $\mathcal{R}_2/I$ is a finite-dimensional (left) vector space over $K(x,y)$, then there exists a nonzero element in the intersection of $I$ and $K(x)\langle D_x, D_y \rangle$.*

Before proving Lemma 1.1, we recall some basic facts about differential operators. Let $\mathcal{A}_2$ be the Weyl algebra $k[x,y]\langle D_x, D_y \rangle$, which is a subring

---
[*]S. Chen was supported by NFS grant CCF-1017217. Z. Li was supported by a grant of the National Natural Science Foundation of China (No. 60821002).
[†]*Emails:* schen21@ncsu.edu (Shaoshi Chen), zmli@mmrc.iss.ac.cn (Ziming Li).

of $\mathcal{R}_2$. Assume that $A$ and $B$ are two nonzero differential operators of the form

$$A = LD_x^m + A_{m-1}D_x^{m-1} + \cdots + A_0 \quad \text{and} \quad B = LD_y^n + B_{n-1}D_x^{n-1} + \cdots + B_0 \tag{1}$$

where $L, A_i, B_j$ are in $K[x, y]$ with $L \neq 0$, and $m, n$ are positive integers. So $A$ and $B$ are in $k[x, y]\langle D_x \rangle$ and $k[x, y]\langle D_y \rangle$, respectively. These two subrings are both contained in $\mathcal{A}_2$.

Leibniz's formula for differentiation is translated into the language of differential operators as: for all $f \in K(x, y)$

$$D_x^k f = \sum_{\ell=0}^{k} \binom{k}{\ell} \frac{\partial^\ell f}{\partial x^\ell} D_x^{k-\ell} \tag{2}$$

and

$$fD_x^k = \sum_{\ell=0}^{k} (-1)^\ell \binom{k}{\ell} D_x^{k-\ell} \frac{\partial^\ell f}{\partial x^\ell}. \tag{3}$$

The relation (2) can be proved by a straightforward induction, while (3) can be proved by applying the adjoint map to (2). Of course, both (2) and (3) hold when $x$ is replaced by $y$. In the sequel, we merely use the facts that, for all $f \in K[x, y]$,

$$D_x^k f = fD_x^k + P \quad \text{and} \quad fD_x^k = D_x^k f - P \tag{4}$$

where $P \in K[x, y]\langle D_x \rangle$ is of degree in $D_x$ less than $k$ and total degree in $x, y$ less than that of $f$.

Let $D = D_x^\beta D_y^\gamma$. If $\beta > m$, Lipshitz claimed that one can always obtain

$$LD \equiv \sum P_\delta D_\delta \mod \langle A \rangle, \tag{5}$$

where the sum on the right hand side is over $D_\delta = D_x^{\delta_1} D_y^{\delta_2}$ with $\delta_1 < \beta$ and $\delta_2 \leq \gamma$. This claim seems not completely correct. In fact, when $\deg_x(L) > 0$ and both $\beta$ and $\gamma$ are positive, multiplying one $L$ is not sufficient to obtain (5). For example, let $D = D_x^m D_y$. Write $A = LD_x^m - R_0$, where $R_0$ is sum of lower order terms in $D_x$. Then

$$LD_x^m \equiv R_0 \mod \langle A \rangle. \tag{6}$$

Multiplying both sides of (6) by $D_y$ yields

$$\begin{align}
D_y L D_x^m &\equiv D_y R_0 \mod \langle A \rangle, \tag{7}\\
LD_y D_x^m - L_y D_x^m &\equiv D_y R_0 \mod \langle A \rangle, \tag{8}\\
LD_y D_x^m &\equiv L_y D_x^m + D_y R_0 \mod \langle A \rangle. \tag{9}
\end{align}$$

In order to reduce the order of $L_y D_x^m$ in (9), we multiply both sides of (9) by $L$, then
$$L^2 D_x^m D_y \equiv \sum P_\delta D_\delta \mod \langle A \rangle,$$
where the sum on the right hand side is over $D_\delta = D_x^{\delta_1} D_y^{\delta_2}$ with $\delta_1 < m$ and $\delta_2 \leq 1$. In contrast to the statement from line $-8$ to $-3$ on page 375 of [1], we have

**Lemma 1.2** *Let $A$ and $B$ be given by (1), and $J$ the left ideal generated by $A$ and $B$ in $\mathcal{A}_2$. Assume that $d$ is an upper bound for the total degrees of $L$, $A_i$ and $B_j$ for all $i,j$ with $0 \leq i \leq m-1$ and $0 \leq j \leq n-1$. Then, for all $\alpha$, $\beta$ in $\mathbb{N}$, we have*
$$L D_x^\alpha D_y^\beta \equiv \sum_{i,j} R_{i,j}^{(\alpha,\beta)} D_x^i D_y^j \mod J,$$
*where $R_{ij} \in K[x,y]$, $\deg R_{i,j}^{(\alpha,\beta)} \leq d$, either $0 \leq i \leq m-1$ and $0 \leq j \leq n-1$ or $i+j \leq \alpha + \beta - 1$.*

**Proof.** If $\alpha < m$ and $\beta < n$, the claim holds. Assume that $\beta \geq n$. We compute
$$\begin{aligned}
L D_x^\alpha D_y^\beta &= (L D_x^\alpha) D_y^\beta = (D_x^\alpha L + P_\alpha) D_y^\beta \quad \text{(by (4))} \\
&= D_x^\alpha L D_y^\beta + P_\alpha D_y^\beta = D_x^\alpha \left( L D_y^{\beta-n} \right) D_y^n + P_\alpha D_y^\beta \\
&= D_x^\alpha \left( D_y^{\beta-n} L + Q_\beta \right) D_y^n + P_\alpha D_y^\beta \quad \text{(by (4))} \\
&= D_x^\alpha D_y^{\beta-n} (L D_y^n) + D_x^\alpha Q_\beta D_y^n + P_\alpha D_y^\beta \\
&= D_x^\alpha D_y^{\beta-n} \left( B - \sum_{j=0}^{n-1} B_j D_y^j \right) + D_x^\alpha Q_\beta D_y^n + P_\alpha D_y^\beta \\
&\equiv -D_x^\alpha D_y^{\beta-n} \left( \sum_{j=0}^{n-1} B_j D_y^j \right) + D_x^\alpha Q_\beta D_y^n + P_\alpha D_y^\beta \mod J.
\end{aligned}$$

It follows from the degree constraints on $P_\alpha$ and $Q_\beta$ that the lemma holds for $\beta \geq n$. Likewise, the lemma holds $\alpha \geq m$. $\square$

Similar to the statement made in line 2 on page 376 in [1], we have

**Lemma 1.3** *Let $A$ and $B$ be given by (1), and $J$ the left ideal generated by $A$ and $B$ in $\mathcal{A}_2$. Assume that $d$ is an upper bound for the total degrees of $L$, $A_i$ and $B_j$ for all $i,j$ with $0 \leq i \leq m-1$ and $0 \leq j \leq n-1$. Then, for all $\alpha$, $\beta$ in $\mathbb{N}$, we have*
$$L^{\alpha+\beta+1-\min(m,n)} D_x^\alpha D_y^\beta \equiv \sum_{i=0}^{m-1} \sum_{j=0}^{n-1} R_{i,j}^{(\alpha,\beta)} D_x^i D_y^j \mod J,$$
*where $R_{ij} \in K[x,y]$ and $\deg R_{i,j}^{(\alpha,\beta)} \leq (\alpha+\beta+1-\min(m,n))\,d$.*

**Proof.** By Lemma 1.2, there are $P_{ij}$, $Q_{ij}$ and $R_{ij}$ in $K[x,y]$ with total degree no more than $d$ such that

$$\begin{aligned} LD_x^\alpha D_y^\beta &\equiv \sum_{i=0}^{m-1} \sum_{j=0}^{n-1} P_{ij} D_x^i D_y^j \\ &+ \sum_{i \geq m, \, 0 \leq i+j \leq \alpha+\beta-1} Q_{ij} D_x^i D_y^j \\ &+ \sum_{j \geq n, \, 0 \leq i+j \leq \alpha+\beta-1} R_{ij} D_x^i D_y^j \quad \mod J. \end{aligned}$$

It follows that

$$\begin{aligned} L^2 D_x^\alpha D_y^\beta &\equiv \sum_{i=0}^{m-1} \sum_{j=0}^{n-1} LP_{ij} D_x^i D_y^j \\ &+ \sum_{i \geq m, \, 0 \leq i+j \leq \alpha+\beta-1} Q_{ij} \left( LD_x^i D_y^j \right) \\ &+ \sum_{j \geq n, \, 0 \leq i+j \leq \alpha+\beta-1} R_{ij} \left( LD_x^i D_y^j \right) \quad \mod J. \end{aligned}$$

Applying Lemma 1.2 to each $LD_x^i D_y^j$ appearing in the second and third summations yields that

$$\begin{aligned} L^2 D_x^\alpha D_y^\beta &\equiv \sum_{i=0}^{m-1} \sum_{j=0}^{n-1} P'_{ij} D_x^i D_y^j \\ &+ \sum_{i \geq m, \, 0 \leq i+j \leq \alpha+\beta-2} Q'_{ij} D_x^i D_y^j \\ &+ \sum_{j \geq n, \, 0 \leq i+j \leq \alpha+\beta-2} R'_{ij} D_x^i D_y^j \quad \mod J \end{aligned}$$

for some $P'_{ij}$, $Q'_{ij}$ and $R'_{ij}$ in $K[x,y]$ with total degrees no more than $2d$. A straightforward induction shows that

$$\begin{aligned} L^k D_x^\alpha D_y^\beta &\equiv \sum_{i=0}^{m-1} \sum_{j=0}^{n-1} P^*_{ij} D_x^i D_y^j \\ &+ \sum_{i \geq m, \, 0 \leq i+j \leq \alpha+\beta-k} Q^*_{ij} D_x^i D_y^j \\ &+ \sum_{j \geq n, \, 0 \leq i+j \leq \alpha+\beta-k} R^*_{ij} D_x^i D_y^j \quad \mod J \end{aligned}$$

for some $P^*_{ij}$, $Q^*_{ij}$ and $R^*_{ij}$ in $K[x,y]$ with total degrees no more than $kd$.

Setting $k = \alpha + \beta + 1 - \min(m,n)$ yields the lemma. $\square$

We are ready to prove Proposition 1.1. Assume further that $I$ is nontrivial. Then $I$ contains two differential operators $A$ and $B$ given by (1). Assume that $J$ is the left ideal generated by $A$ and $B$ in $\mathcal{A}_2$, It suffices to show that there is a nonzero element in the intersection of $J$ and $K[x]\langle D_x, D_y \rangle$.

We apply the same counting argument used in [1]. Assume that $d$ is an upper bound for all coefficients $A_i$ and $B_j$ Let $N$ a positive integer, and let

$$V_N = \left\{ L^N x^\gamma D_x^\alpha D_y^\beta \mid \gamma, \alpha, \beta \in \mathbb{N}, \, \gamma + \alpha + \beta \leq N \right\}$$

and

$$W_N = \left\{ x^s y^t D_x^i D_y^j \mid s, t, i, j \in \mathbb{N}, \, s+t \leq N(d+1), \, i < m, \, j < n \right\}.$$

By Lemma 1.3, $L^N x^\gamma D_x^\alpha D_y^\beta$ is congruent a $K$-linear combination of the elements in $W_N$ modulo $J$. Since $|V_N| = O(N^3)$ and $|W_N| = O(N^2)$. there must be a nontrivial $K$-linear combination of the elements in $V_N$ lying in $J$ when $N$ is sufficiently large. □